\documentclass[12pt]{article} 
\usepackage{amssymb}
\def\eps{\varepsilon} \def\fin{_{\rm fin}} \def\opp{^{\rm opp}}
\def\id{{\rm id}} \def\cir{\,{\scriptstyle \circ}\,}
\def\eins{\hbox{\rm 1}\mskip-4.4mu\hbox{\rm l}}
\def\LI#1#2{\Phi^#2_{#1_#2}} \def\RI#1#2{\Psi^#2_{#1_#2}} 
\def\A#1#2{{\alpha^#2_{#1_#2}}} 
\def\CC{{\mathbb C}} \def\TT{{\cal T}} \def\RR{{\mathbb R}} 
\def\NN{{\mathbb N}} 
\def\QED{\hspace*{\fill}$\square$} 
\newcounter{nummer}[section]
\newenvironment{satz}[1]{\addtocounter{nummer}{1}\vskip5pt
       \noindent\bf #1 \arabic{section}.\arabic{nummer}. \it}{\rm\vskip3pt}
\newenvironment{beweis}[1]{\vskip3pt\noindent\it Proof#1.\rm}{\QED\vskip3pt}
\newenvironment{acknowledgement}{\vskip7mm
\noindent{\large\bf Acknowledgement}\vskip3mm\rm}{\vskip3mm}
\begin{document}
\author{K.-H. Rehren \\[1mm] 
{\small Institut f\"ur Theoretische Physik, Universit\"at 
G\"ottingen} \\ {\small Bunsenstra{\ss}e 9, D-37073 G\"ottingen} \\ 
{\small E-mail: rehren@theorie.physik.uni-goettingen.de}}
\date{} 
\title{Canonical Tensor Product Subfactors}

\maketitle 
\begin{abstract}
Canonical tensor product subfactors (CTPS's) describe, among other things, the
embedding of chiral observables in two-dimensional conformal quantum
field theories. A new class of CTPS's is constructed some of which 
are associated with certain modular invariants, thereby
establishing the expected existence of the corresponding
two-dimensional theories.
\end{abstract}
\def\Large{\large}

\section{Introduction and results}

There is a common mathematical note which recurs again and again in
the areas of conformal quantum field theory and modular invariants on
the one hand, and asymptotic subfactors and quantum doubles on the
other hand. In all these areas, there arise inclusions of von Neumann
factors of the form $A\otimes B\subset C$ sharing a ``canonical''
property (see Def.\ 1.1 below) for which we call them ``canonical
tensor product subfactors'' (CTPS, cf.\ \cite{R1}). E.g., in chiral quantum
field theories on $S^1$, CTPS's describe the violation of Haag duality
for disjoint intervals (Jones-Wassermann subfactors, cf.\ \cite{KLM,X2}),
or the embedding of a coset model into a given ambient model. In
two-dimensional conformal quantum field theories they describe the
embedding of chiral subtheories \cite{R1} which is (incompletely)
reflected also by modular invariant coupling matrices \cite{MS}. 
Ocneanu's asymptotic subfactors \cite{O} which are sometimes regarded as
generalized quantum doubles \cite{EK} are also CTPS's. 

The main result in this article is the presentation (Thm.\ 1.4) of a class
of new CTPS's associated with extensions of closed systems
of endomorphisms (Def.\ 1.2, 1.3). Among them there is a 
subclass of considerable importance for the understanding of modular 
invariants. Namely, with every modular invariant constructed by a method 
due to B\"ockenhauer, Evans and Kawahigashi \cite{BEK} one can
associate one of the new CTPS's which, if 
interpreted as a local inclusion of an algebra of chiral observables
into an algebra of two-dimensional observables, allows to prove the 
existence of a complete two-dimensional local conformal quantum field
theory associated with the given modular invariant (Cor.\ 1.6).  

The mathematical abstraction of this physical problem as a
problem on von Neumann algebras and subfactors is most efficient. It
is based on the seminal realization \cite{H} that positive-energy
representations (``superselection charges'') and their fusion are
conveniently expressed in terms of endomorphisms, promoting particle
statistics to a unitary operator representation (braiding) on the
physical Hilbert space, and identifying the statistical dimension as
(the square root of) a Jones index. 

Extensions or embeddings of quantum field theories can also be coded
into single subfactors \cite{LR}. The characterization of a
subfactor, in turn, in terms of a ``Q-system'' \cite{L} is
particularly useful in this context since these data directly describe
the charged field content of the extended theory in terms of
superselection charges of the embedded theory \cite{LR,RST}. The new
subfactors presented in Thm.\ 1.4 are also defined by specification of
their Q-systems (in terms of certain matrix elements for the
transition between two extensions), thus making as close contact with
the structure of modular invariants as possible.   

CTPS's are very special cases of ``symmetric joint inclusions'', i.e.,
triples of von Neumann algebras $(A,B,C)$ such that $A$ and $B$ are
commuting subalgebras of $C$. After a survey of some general properties of
symmetric joint inclusions in Sect.\ 4, we give a characterization of
``normality'' for CTPS's in Proposition 4.3. This is a maximality
property which, in the case of the embedding of chiral observables into
a two-dimensional conformal quantum field theory, corresponds to the
maximally extended chiral algebras and diagonal or permutation
invariants \cite{R1}. 

The canonical property mentioned before is a natural feature
of the embedding of chiral quantum field theories into a
two-dimensional conformal quantum field theory, reflecting the
independence of left and right moving degrees of freedom \cite{R1}. It
is defined as follows. 
\begin{satz}{Definition} 
  A tensor product subfactor of the form $A\otimes B \subset C$ is
  called a {\em canonical tensor product subfactor (CTPS)} 
  if either $A,B,C$ are type II factors and $C$ considered as an 
  $A\otimes B$-$A\otimes B$ bimodule decomposes into irreducibles
  which are all tensor products of $A$-$A$ bimodules with $B$-$B$ bimodules,
  or if $A,B,C$ are type III factors and the dual canonical endomorphism
  $\theta\equiv\bar\iota\cir\iota\in End(A\otimes B)$ decomposes into
  irreducibles which are all tensor products of endomorphisms of $A$ with 
  endomorphisms of $B$. Let, in the type III case,
  $$ \theta\simeq\bigoplus_{\alpha,\beta}\; Z_{\alpha,\beta} 
  \;\alpha \otimes \beta $$
  where the sum extends over two sets of mutually inequivalent irreducible 
  endomorphisms of $A$ and of $B$, respectively. Then we call the matrix
  of multiplicities $Z$ with non-negative integer entries the 
  {\em coupling matrix} of the CTPS. The coupling matrix in the type II case
  is defined analogously in terms of mutually inequivalent irreducible
  $A$-$A$ and $B$-$B$ bimodules.
\end{satz}
Here, as always throughout this paper, $\iota:A\otimes B\to C$ denotes
the inclusion homomorphism of the subfactor under consideration, and
$\bar\iota:C\to A\otimes B$ a conjugate homomorphism \cite{LRo}.

In order to state our main result, we have to introduce some further 
notions. We consider type III von Neumann factors $N$, and denote by
$End\fin(N)$ the set of unital endomorphisms $\lambda$ of $N$ with
finite dimension $d(\lambda)$. 
\begin{satz}{Definition}
  A {\em closed $N$-system} is a set $\Delta \subset End\fin(N)$ of mutually 
  inequivalent irreducible endomorphisms such that (i)
  $\id_N\in\Delta$, (ii) if $\lambda\in\Delta$ then there is a
  conjugate endomorphism $\bar\lambda\in\Delta$, and (iii) if
  $\lambda,\mu\in\Delta$ then $\lambda\mu$ belongs to
  $\Sigma(\Delta)$, the set of endomorphisms which are 
  equivalent to finite direct sums of elements from $\Delta$.
\end{satz}
\begin{satz}{Definition} 
  Let $N\subset M$ be a subfactor with inclusion homomorphism 
  $\iota:N\to M$. An {\em extension of the closed $N$-system $\Delta$} 
  is a pair $(\iota,\alpha)$, where $\iota$ is as above, and $\alpha$ is
  a map $\Delta\to End\fin(M)$, $\lambda \mapsto\alpha_\lambda$, such that 
  \\[2pt]\indent
  (E1) $\quad\iota\cir\lambda = \alpha_\lambda\cir\iota$, \\[2pt]\indent
  (E2) $\quad\iota(Hom(\nu,\lambda\mu))\subset 
  Hom(\alpha_\nu,\alpha_\lambda\alpha_\mu)$.
\end{satz}
After these preliminaries, we can state our main result. 
\begin{satz}{Theorem} 
  Let $N_1\subset M$ and $N_2\subset M$ be two subfactors 
  of $M$, and $(\iota_1,\alpha^1)$ and $(\iota_2,\alpha^2)$ a pair of 
  extensions of a finite closed $N_1$-system $\Delta_1$ 
  and a finite closed $N_2$-system $\Delta_2$, respectively. Then there 
  exists an irreducible CTPS 
  $$ A\equiv N_1\otimes N_2\opp \subset B$$
  with dual canonical endomorphism
  $$ \theta\equiv\bar\iota\cir\iota \simeq
  \bigoplus_{\lambda_1\in\Delta_1,\lambda_2\in\Delta_2}
  Z_{\lambda_1,\lambda_2}\; \lambda_1\otimes\lambda_2\opp , $$
  whose coupling matrix $Z$ of multiplicities is given by 
  $$ Z_{\lambda_1,\lambda_2} = \dim Hom(\A\lambda1,\A\lambda2) . $$
\end{satz}
The special case when $\Delta_i$ are braided systems is of particular
interest for the problem of modular invariants in conformal quantum
field theory:  
\begin{satz}{Proposition} 
  Assume in addition that the closed systems $\Delta_1$ 
  and $\Delta_2$ are braided with unitary braidings $\eps_1$ and
  $\eps_2$, respectively, turning $\Pi(\Delta_1)$ and $\Pi(\Delta_2)$
  into braided monoidal categories. If for all
  $\lambda_i,\mu_i\in\Delta_i$ and all
  $\phi\in Hom(\A\lambda1,\A\lambda2)$, $\psi\in Hom(\A\mu1,\A\mu2)$,
  one has \\[2pt]\indent
  (E3) $\quad(\psi\times\phi)\cir \iota_1(\eps_1(\lambda_1,\mu_1)) = 
  \iota_2(\eps_2(\lambda_2,\mu_2))\cir (\phi\times\psi$), \\[2pt]
  then the canonical isometry $w_1\in Hom(\theta,\theta^2)$ 
  (defined below in the proof of the Theorem) and the braiding operator
  $\eps(\theta,\theta)$ naturally induced by the braidings $\eps_1$ 
  and $\eps_2\opp$ satisfy
  $$ \eps(\theta,\theta) w_1 = w_1. $$
\end{satz}
This result answers an open question in quantum field theory, where possible 
matrices $Z$ are classified which are supposed to describe the 
restriction of a given two-dimensional modular invariant conformal 
quantum field theory to its chiral subtheories, while it is actually
not clear whether any given solution $Z$ does come from a two-dimensional 
quantum field theory. 
This turns out to be true for a large class of solutions, obtained
in \cite{BEK}:
\begin{satz}{Corollary}
  Let $A: I \mapsto A(I)$ be a chiral net of local observables associated
  with the open intervals $I\subset \RR$, such that each $A(I)$ is a type III 
  factor. Let $\Delta_{\rm DHR}$ be a closed system of mutually 
  inequivalent irreducible DHR-endomorphisms \cite{H} of $A$ with
  finite dimension, localized in some interval $I_0$, and put
  $N:=A(I_0)$. We assume that $A$ is conformally covariant, implying
  \cite{GL} that the system of restrictions 
  $\Delta:=\{\lambda=\lambda_{\rm DHR}\vert_N : 
  \lambda_{\rm DHR}\in\Delta_{\rm DHR}\}$ is a closed $N$-system.
  Let $N_1\subset N$ be a subfactor with canonical endomorphism $\theta
  \in \Sigma(\Delta)$, and $N\subset M$ its Jones extension. Put 
$$Z_{\lambda,\mu}:= \dim Hom(\alpha^+_\lambda,\alpha^-_\mu)$$
  where $\alpha^\pm$ are the pair of $\alpha$-inductions \cite{LR,BE,X1} 
  of endomorphisms of $N$ to endomorphisms of $M$, associated with 
  the braidings given by the DHR statistics and its opposite. Then 
  there is a two-dimensional local conformal quantum field theory
  described by a net $B:O\mapsto B(O)$ of observables associated
  with the double-cones $O=I\times J$ in $\RR^2$, containing subnets
  of left and right chiral observables $A_L$ and $A_R$ both isomorphic
  with $A$, such that the local inclusions of chiral observables
  $A_L(I)\times A_R(J)\subset B(O)$ are CTPS's with coupling matrix
  $ZC$. (Here $C$ is the matrix describing sector conjugation in
  $\Delta$.) Equivalently, the restriction of the vacuum
  representation of the two-dimensional quantum field theory $B$ to
  its chiral subtheories is given by 
$$\pi^0\vert_{A_L\otimes A_R} = \sum_{\lambda,\mu\in\Delta}\;
  Z_{\lambda,\bar\mu}\;\pi_\lambda\otimes\pi_\mu.$$
\end{satz}
The corollary combines and adapts results from \cite{BEK,LR}. The
point is that if the dual canonical endomorphism $\theta$ associated
with $N\subset M$ belongs to $\Sigma(\Delta)$, then $\alpha$-induction
\cite{LR,BE,X1} provides a pair of extensions $(\iota,\alpha^+)$ and
$(\iota,\alpha^-)$ which satisfies (E1), (E2) as well as (E3) (e.g., 
\cite[I; Def.\ 3.3, Lemma 3.5 and 3.25]{BE}). The associated coupling matrix 
$Z_{\lambda,\mu}$ is automatically a modular invariant \cite{BEK}.  
By the characterization of extensions of local quantum field theories 
given in \cite[Prop.\ 4.9]{LR}, the local subfactor 
$A_L(I_0)\otimes A_R(I_0) \subset B(O_0)$ given by the Thm.\ 1.4
induces an entire net of subfactors, indexed by the double-cones $O$
of two-dimensional Minkowski space. (The charge conjugation $C$ arises
due to an anti-isomorphism between $N\opp$ and $N$, cf.\ 
\cite[Prop.\ 4.10{\it ff}\/]{LR}.) The statement of Proposition 1.5 is
precisely the criterium given in \cite{LR} for the resulting
two-dimensional quantum field theory to be local.

Thus, every modular invariant found by the $\alpha$-induction method given in
\cite{BEK} indeed corresponds to a two-dimensional local conformal
quantum field theory extending the given chiral nets of observables. 

\section{Extensions of systems of endomorphisms}
We collect some immediate consequences of the definition  
of an extension, Def.\ 1.3, using terminology and notations as in
\cite{DR,LRo}. 
\begin{satz}{Proposition}
  An extension $(\iota,\alpha)$ of a closed $N$-system gives rise to
  a monoidal functor from
  the full monoidal C* subcategory of $End\fin(N)$ with objects 
  $\Pi(\Delta)$ (the set of finite products of elements from $\Delta$)
  into the monoidal C* category $End\fin(M)$. In particular, 
  $\lambda\mu\simeq \bigoplus_\nu N_{\lambda\mu}^\nu\,\nu$ implies
  $\alpha_\lambda\alpha_\mu\simeq \bigoplus_\nu
  N_{\lambda\mu}^\nu\,\alpha_\nu$ (notwithstanding $\alpha_\lambda$ will
  be reducible in general), and $\alpha_{\bar\lambda}$ is conjugate
  to $\alpha_\lambda$. 
\end{satz}
\begin{beweis}{}
The functor maps objects $\lambda_1\cir\dots\cir\lambda_n$ to
$\alpha_{\lambda_1}\cir\dots\cir\alpha_{\lambda_n}$, and 
intertwiners $T$ to $\iota(T)$ which are again intertwiners by iteration 
of (E2). It follows from (E2) that $\alpha$ preserves the fusion rules
as stated.
In particular, $\alpha_{\id_N}$ is an idempotent within $End\fin(M)$,
implying that its dimension is $1$, hence it is invertible and must
coincide with $\id_M$. Thus the functor preserves the monoidal 
unit object. It preserves the right monoidal product of intertwiners
trivially, and the left monoidal product by (E1). 
Conjugacy between $\alpha_{\bar\lambda}$ and $\alpha_\lambda$ 
is a consequence of the following lemma.
\end{beweis}
\begin{satz}{Lemma} Let   $(\iota,\alpha)$ be an extension of a closed 
  $N$-system $\Delta$. Let $R\in Hom(\id_N,\bar\lambda\lambda)$
  and $\bar R\in Hom(\id_N,\lambda\bar\lambda)$ be a pair of standard
  isometries solving the conjugate equations
  $$(1_\lambda\times R_\lambda^*)\cir
  (\bar R_\lambda\times 1_\lambda)=d(\lambda)^{-1}1_\lambda
  = (1_{\bar\lambda}\times \bar R_\lambda^*)\cir 
  (R_\lambda\times 1_{\bar\lambda}),$$ 
  and thus implementing the unique left- and right-inverses \cite{LRo} 
  $\Phi_\lambda$ and $\Psi_\lambda$ for $\lambda\in\Delta$. 
  Then $\iota(R_\lambda)$ and $\iota(\bar R_\lambda)$ induce left- and 
  right-inverses $\Phi_{\alpha_\lambda}$ and $\Psi_{\alpha_\lambda}$ for
  $\alpha_\lambda$. If either $N\subset M$ has finite index, or
  $\Delta$ is a finite system, then $d(\alpha_\lambda)=d(\lambda)$,
  and $\Phi_{\alpha_\lambda}$ and $\Psi_{\alpha_\lambda}$ are the unique
  standard left- and right-inverses.
\end{satz}
\begin{beweis}{}
The first statement is an obvious consequence of (E2). If the index
$d(\iota)^2$ is finite, then (E1) implies
$d(\alpha_\lambda)=d(\lambda)$. If $\Delta$ is finite, then the minimal 
dimensions $d(\alpha_\lambda)$ are uniquely determined by the fusion 
rules of $\{\alpha_\lambda,\lambda\in\Delta\}$, and the latter
coincide with those of $\{\lambda\in\Delta\}$. Hence again
$d(\alpha_\lambda)=d(\lambda)$. Since $d(\lambda)$ are also the
dimensions associated with the pair of isometries $\iota(R_\lambda)$,
$\iota(\bar R_\lambda)$, the last claim follows by \cite[Thm.\ 3.11]{LRo}.  
\end{beweis}
Thus, general properties of standard left- and right-inverses \cite{LRo} 
are applicable. We shall in the sequel repeatedly exploit the trace
property 
$$d(\rho)\Phi_\rho(S^*T)=d(\tau)\Phi_\tau(TS^*)
\qquad {\rm if}\qquad S,T\in Hom(\rho,\tau) $$
for standard left-inverses of $\rho,\tau\in End\fin(M)$, their
multiplicativity $\Phi_{\tau\rho}=\Phi_\tau\Phi_\rho$, 
as well as the equality of standard left- and right-inverses
$\Phi_\rho=\Psi_\rho$ on $Hom(\rho,\rho)$. 
\section{Construction of the new CTPS's}
We shall prove Theorem 1.4 by the specification of ``Q-systems'' (or
``canonical triples'') $(\theta,w,w_1)$, which uniquely determine
subfactors \cite{L}. 

Longo's characterization states that $\theta\in End\fin(A)$ is the dual 
canonical endomorphism associated with $A\subset B$ if (and only if)
there is a pair of isometries $w\in Hom(\id_A,\theta)$ and 
$w_1\in Hom(\theta,\theta^2)$ satisfying\\[2pt]
\indent
(Q1) $\quad w^*w_1 = \theta(w^*)w_1 = d(\theta)^{-1/2} \eins_A$, 

(Q2) $\quad w_1w_1 = \theta(w_1)w_1$, and 

(Q3) $\quad w_1w_1^* = \theta(w_1^*)w_1$.\\[2pt]
Namely, then the map $w_1^*\theta(\,\cdot\,)w_1$ is the minimal
conditional expectation onto its image $A_1=w_1^*\theta(A)w_1\subset
A$. For $\iota_1:A_1\to A$ the inclusion map and $\bar\iota_1: A\to A_1$ 
defined by $\theta=\iota_1\bar\iota_1$, the pair of isometries
$w\in Hom(\id_A,\iota_1\bar\iota_1)$ and $\iota_1^{-1}(w_1)\in
Hom(\id_{A_1},\bar\iota_1\iota_1)$ achieves the conjugacy between
$\iota_1$ and $\bar\iota_1$. By the Jones construction \cite{J}, then,
the subfactor $A_1\subset A$ determines its dual subfactor (the Jones
extension) $A\subset B$ such that $\theta=\bar\iota\iota$.  

\begin{beweis}{ of Theorem 1.4}
First notice that the multiplicity of $\id_A$ in $\theta$ is
$Z_{\id_{N_1},\id_{N_2}} = \dim Hom(\id_M,\id_M) = 1$, so the asserted
subfactor is automatically irreducible. 

In order to show that $\theta$ given in the Theorem is the dual
canonical endomorphism associated with a subfactor $A\subset B$, we
construct the Q-system $(\theta,w,w_1)$ as follows.
We first choose a complete system of mutually inequivalent isometries  
$W_{(\lambda_1,\lambda_2,l)}\equiv W_l\in A\equiv N\otimes N\opp$, 
where $l$ is considered as a multi-index
$(\lambda_1\in\Delta_1,\lambda_2\in\Delta_2,
l=1,\dots Z_{\lambda_1,\lambda_2})$, and put 
$$\theta=\sum_l W_l\; (\lambda_1\otimes \lambda_2\opp)(\,\cdot\,)\; W_l^*. $$
The choice of these isometries is immaterial and affects the subfactor
to be constructed only by inner conjugation.

Since $Hom(\id_A,\theta)$ is one-dimensional, the isometry $w$ is already
fixed up to an irrelevant complex phase: $w=W_0$, where $0$
refers to the multi-index $l=0\equiv(\id_{N_1},\id_{N_2},1)$.
The second isometry, $w_1$, must be of the form
$$ w_1=\sum_{l,m,n} (W_l \times W_m) \cir \TT_{lm}^n \cir W_n^* $$
where $\TT_{lm}^n\in Hom(\nu_1\otimes\nu_2\opp,
(\lambda_1\otimes\lambda_2\opp)\cir(\mu_1\otimes\mu_2\opp))$, since these
operators span $Hom(\theta,\theta^2)$. 

In turn, $\TT_{lm}^n$ must be of the form 
$$ \TT_{lm}^n = \sum_{e_1,e_2} \zeta_{lm,e_1e_2}^n \; T_{e_1} \otimes
(T_{e_2}^*)\opp \qquad (\zeta_{lm,e_1e_2}^n\in\CC) $$
where $T_{e_i}$ constitute orthonormal bases of the intertwiner
spaces $Hom(\nu_i,\lambda_i\mu_i)$, since these operators span
the spaces $Hom(\nu_1\otimes\nu_2\opp,
(\lambda_1\otimes\lambda_2\opp)\cir(\mu_1\otimes\mu_2\opp))\equiv
Hom(\nu_1,\lambda_1\mu_1)\otimes Hom(\nu_2\opp,\lambda_2\opp\mu_2\opp)$. 
Note that if $T\in Hom(\alpha,\beta)$ is isometric in $N$, then
$(T^*)\opp \in Hom(\beta,\alpha)\opp \equiv Hom(\alpha\opp,\beta\opp)$ is
isometric in $N\opp$. The labels $e_i$ are again multi-indices
$(\lambda,\mu,\nu,e=1,\dots \dim Hom(\nu,\lambda\mu))$. 

It remains therefore to determine the complex coefficients 
$\zeta_{lm,e_1e_2}^n$, such that $w_1$ is an isometry satisfying Longo's
relations (Q1--3) above. To specify these coefficients, we equip the
spaces $Hom(\A\lambda1,\A\lambda2)$ with the non-degenerate scalar products 
$(\phi,\phi'):=\LI\lambda1(\phi^*\phi')$ (where $\LI\lambda i$ stand
for the induced left-inverses for $\A\lambda i$, $i=1,2$, cf.\ Lemma 2.2). 
With respect to these scalar products, we choose orthonormal bases
$\{\phi_l,l=1,\dots Z_{\lambda_1,\lambda_2}\}$ for all
$\lambda_1\in\Delta_1$, $\lambda_2\in\Delta_2$, and put 
$$ \zeta_{lm,e_1e_2}^n = \sqrt\frac{d(\lambda_2)d(\mu_2)}{d(\theta)d(\nu_2)}
\; \LI\nu1
[\iota_1(T_{e_1}^*)(\phi_l^* \times \phi_m^*)\iota_2(T_{e_2})\phi_n]. $$
This formula is only apparently asymmetric under exchange 
$1\leftrightarrow 2$: by the trace property $d(\lambda_2)
\LI\lambda2(\phi\phi^*) = d(\lambda_1)\LI\lambda1
(\phi^*\phi)$, an orthonormal basis $\psi_l$ of
$Hom(\A\lambda2,\A\lambda1)$ differs from $\phi_l^*$ by a factor
$\sqrt\frac{d(\lambda_1)}{d(\lambda_2)}$, so that in fact  
$$\overline{\zeta_{lm,e_1e_2}^n} = 
\sqrt\frac{d(\lambda_1)d(\mu_1)}{d(\theta)d(\nu_1)}\;
\LI\nu2
[\iota_2(T_{e_2}^*)(\psi_l^* \times \psi_m^*)\iota_1(T_{e_1})\psi_n] . $$

With these coefficients, condition (Q1) is trivially satisfied, since
left multiplication of $w_1$ by $w^*$ singles out the term $l=0$ due
to $W_0^*W_l=\delta_{l0}$. This leaves only terms with
$\lambda_i=\id_{N_i}$, hence $\mu_i=\nu_i$, for which  $T_{e_i}$ are
trivial and $\sqrt{d(\theta)}\zeta_{0m,e_1e_2}^n= \delta_{mn}$ (up to
cancelling complex phases), so  
$\sqrt{d(\theta)}w^*w_1 = \sum_n W_nW_n^*=\eins_A$. For
$\theta(w^*)w_1$ the argument is essentially the same.

We turn to the conditions (Q2) and (Q3). Whenever we compute either of
the four products occurring, we obtain a Kronecker delta
$W_s^*W_t=\delta_{st}$ for one pair of the labels $l,m,n,\dots$
involved, while the remaining operator parts are of the form 
\begin{eqnarray*}
(W_l \times W_m \times W_k) \left[(T_{e_1} \times 1_{\kappa_1})T_{f_1}
\otimes (((T_{e_2} \times 1_{\kappa_2})T_{f_2})^*)\opp\right] W_n^*, \\[2pt] 
(W_l \times W_m \times W_k) \left[(1_{\lambda_1} \times T_{g_1})T_{h_1}
\otimes (((1_{\lambda_2} \times T_{g_2})T_{h_2})^*)\opp\right] W_n^* \;
\end{eqnarray*}
for the left- and right-hand side of (Q2), $w_1w_1 = \theta(w_1)w_1$,
and in turn, 
\begin{eqnarray*} 
(W_l \times W_m) \left[T_{e_1}T_{f_1}^* \otimes ((T_{e_2}T_{f_2}^*)^*)\opp\right] 
(W_n \times W_k)^* , \qquad\qquad\qquad \\[2pt]
(W_l \times W_m) \left[(1_{\lambda_1}\!\times\!
  T_{g_1}^*)(T_{h_1}\!\times\! 1_{\kappa_1}) \otimes 
(((1_{\lambda_2}\!\times\! T_{g_2}^*)(T_{h_2}\!\times\!
1_{\kappa_2}))^*)\opp\right] 
(W_n \times W_k)^*  \end{eqnarray*}
for the left- and right-hand side of (Q3), $w_1w_1^* = \theta(w_1^*)w_1$.
(In these expressions, we do not specify the respective intertwiner
spaces to which the various operators $T$ belong, since 
these are determined by the context.)

The numerical coefficients of these operators are sums over products of two
$\zeta$'s or one $\zeta$ and one $\overline{\zeta}$, respectively,
with a summation over one common label $s=1,\dots Z_{\sigma_1,\sigma_2}$ 
due to the above Kronecker $\delta_{st}$. These sums can be carried out.

Namely, the
coefficients of the above operators on both sides of (Q2) involve one
factor $\zeta^s_{\dots}$ which is a scalar product of the form 
$\LI\sigma1(X\phi_s)=(X^*,\phi_s)$ in $Hom(\A\sigma1,\A\sigma2)$, 
so summation with the operator $\phi_s^*$ contributing to the other 
factor $\zeta$ yields $\sum_s$ $(X^*,\phi_s)\phi_s^* = X$. Then the 
coefficients of the above operators on both sides of (Q2) are easily cast
into the respective form 
\begin{eqnarray*}
\sum_s\zeta^s_{lm,e_1e_2}\zeta^n_{sk,f_1f_2}
=\sqrt\frac{d(\lambda_2)d(\mu_2)d(\kappa_2)}{d(\theta)^2d(\nu_2)}\times
\qquad\qquad\qquad\qquad\qquad\qquad \\
\LI\nu1[\iota_1(T_{f_1}^*(T_{e_1}^*\times 1_{\kappa_1}))
(\phi_l^*\times\phi_m^*\times\phi_k^*)
\iota_2((T_{e_2}\times 1_{\kappa_2})T_{f_2})
\phi_n], \\
\sum_s\zeta^s_{mk,g_1g_2}\zeta^n_{ls,h_1h_2}=
\sqrt\frac{d(\lambda_2)d(\mu_2)d(\kappa_2)}{d(\theta)^2d(\nu_2)}\times
\qquad\qquad\qquad\qquad\qquad\qquad \\
\LI\nu1[\iota_1(T_{h_1}^*(1_{\lambda_1}\times T_{g_1}^*))
(\phi_l^*\times\phi_m^*\times\phi_k^*)
\iota_2((1_{\lambda_2}\times T_{g_2})T_{h_2})
\phi_n].\end{eqnarray*}
Now, since the passage from bases of the form $(T_e\times 1_\kappa)T_f$
to bases $(1_\lambda\times T_g)T_h$ of $Hom(\nu,\lambda\mu\kappa)$ for
any fixed $\nu,\lambda,\mu,\kappa$ is described by unitary matrices,
equality of both sides of (Q2) follows. 

The case of (Q3) is in the same vein, but slightly more involved. In the
coefficients on the left-hand 
$\sum_s\zeta^s_{lm,e_1e_2}\overline{\zeta^s_{nk,f_1f_2}}$, 
we read again the first factor as a scalar product $(X^*,\phi_s)$
within $Hom(\A\sigma1,\A\sigma2)$ and perform the summation $\sum_s
(X^*,\phi_s)\phi_s^* = X$ with the operator $\phi_s^*$ contributing to
the second factor. 
This yields, after application of the trace property for standard 
left-inverses, the coefficients on the left-hand side of (Q3)
\begin{eqnarray*}\sum_s\zeta^s_{lm,e_1e_2}\overline{\zeta^s_{nk,f_1f_2}}=
\sqrt\frac{d(\lambda_2)d(\mu_2)d(\kappa_2)d(\nu_2)}
{d(\theta)^2d(\sigma_2)^2}
\frac{d(\lambda_1)d(\mu_1)}{d(\sigma_1)}\times \quad\qquad\qquad \\
\LI{\mu_1\lambda}1[(\phi_l^*\times\phi_m^*)
\iota_2(T_{e_2}T_{f_2}^*)
(\phi_n\times\phi_k)
\iota_1(T_{f_1}T_{e_1}^*)].\end{eqnarray*}

To compute the coefficients 
$\sum_s\overline{\zeta^m_{sk,g_1g_2}}\zeta^n_{ls,h_1h_2}$ on the
right-hand side of (Q3), we first rewrite the second factor
as a scalar product $(\phi_s,X)$ within $Hom(\A\sigma1,\A\sigma2)$. 
This is achieved by applying the trace property:
$$\zeta^n_{ls,h_1h_2} = \sqrt\frac{d(\lambda_2)d(\sigma_2)}{d(\theta)d(\nu_2)}
\frac{d(\lambda_1)d(\sigma_1)}{d(\nu_1)}\;
\LI\sigma1[\phi_s^*\LI\lambda1((\phi_l^*\times 1_{\A\sigma2})
\iota_2(T_{h_2})\phi_n\iota_1(T_{h_1}^*))]$$
Now the sum over $s$ with $\phi_s$ contributing to
$\overline{\zeta^m_{sk,g_1g_2}}$ can be performed as before,
yielding the coefficients on the right-hand side of (Q3) in the form
\begin{eqnarray*}
\sum_s\overline{\zeta^m_{sk,g_1g_2}}\zeta^n_{ls,h_1h_2}= 
\sqrt\frac{d(\lambda_2)d(\kappa_2)d(\sigma_2)^2}{d(\theta)^2d(\nu_2)d(\mu_2)}
\frac{d(\lambda_1)d(\sigma_1)}{d(\nu_1)} \times 
\qquad\qquad\qquad\qquad\qquad \\
\LI{\mu_1\lambda}1
[(\phi_l^*\!\times\!\phi_m^*)
\iota_2((1_{\lambda_2}\!\times\! T_{g_2}^*)(T_{h_2}\!\times\! 1_{\kappa_2}))
(\phi_n\!\times\!\phi_k)
\iota_1((T_{h_1}^*\!\times\! 1_{\kappa_1})(1_{\lambda_1}\!\times\! T_{g_1}))].
\end{eqnarray*}
Noting that the passage from bases $\sqrt\frac{d(\mu)}{d(\sigma)}T_eT_f^*$ 
to bases $\sqrt\frac{d(\sigma)}{d(\nu)}(1_{\lambda}\times
T_{g}^*)(T_h\times 1_{\kappa})$ of $Hom(\nu\kappa,\lambda\mu)$ 
is again described by unitary matrices for any fixed
$\nu,\kappa,\lambda,\mu$, we obtain equality of both sides of (Q3).

It remains to show that $w_1$ is an isometry, $w_1^*w_1=1$.

Performing the multiplication $w_1^*w_1$ yields two Kronecker delta's
from the factors $W_l \times W_m$, and two more Kronecker delta's from
the factors $T_{e_1}\otimes (T_{e_2}^*)\opp$. Thus
$$w_1^*w_1=\sum_{ns}\left(\sum_{lm,e_1e_2}
\overline{\zeta^s_{lm,e_1e_2}}\zeta^n_{lm,e_1e_2}\right) W_sW_n^*, $$
and we have to perform the sums over $l,m,e_1,e_2$ (involving, as
sums over multi-indices, the summation over sectors
$\lambda_i,\mu_i\in\Delta_i$ for fixed $\nu_i\in\Delta_i$, $i=1,2$). 

It turns out convenient to express $\zeta^n_{lm,e_1e_2}$ as
a scalar product $(\phi_m,X)$ within $Hom(\A\mu1,\A\mu2)$ as
before (with indices relabelled), and to perform the sum over $m$
first. This yields 
\begin{eqnarray*}\sum_{lm,e_1e_2} 
\overline{\zeta^s_{lm,e_1e_2}}\zeta^n_{lm,e_1e_2}
=\sum_{l,e_1e_2}
\frac{d(\lambda_2)d(\mu_2)}{d(\theta)d(\nu_2)}
\frac{d(\lambda_1)d(\mu_1)}{d(\nu_1)}  \times\qquad\qquad\qquad\qquad \\
\LI\nu1[\phi_s^*\iota_2(T_{e_2}^*)
\left(\phi_l\times\LI\lambda1[(\phi_l^*\times 1_{\A\mu1})
\iota_2(T_{e_2})\phi_n\iota_1(T_{e_1}^*)]\right)\iota_1(T_{e_1})]. 
\end{eqnarray*}
In this expression, we can perform the sums over $(e_1,\mu_1)$ and over 
$(e_2,\mu_2)$ after a unitary passage from the bases of orthonormal
isometries $T_e$ of $Hom(\nu,\lambda\mu)$ to the bases 
$\sqrt\frac{d(\lambda)d(\nu)}{d(\mu)}
(1_\lambda\times T_{e'_1}^*)(\bar R_\lambda\times 1_\nu)$, and obtain
after use of the conjugate equations for $R_\lambda,\bar R_\lambda$ 
$$\sum_{lm,e_1e_2} \overline{\zeta^s_{lm,e_1e_2}}\zeta^n_{lm,e_1e_2} =
\sum_{l,\lambda_1\lambda_2} \frac{d(\lambda_2)^2}{d(\theta)}\; 
\LI\nu1[\RI\lambda2(\phi_l\phi_l^*)\times(\phi_s^*\phi_n)].$$
Here $\RI\lambda2$ ist the standard right-inverse implemented by
$\iota_2(\bar R_{\lambda_2})$ which coincides with $\LI\lambda2$
on $Hom(\A\lambda2,\A\lambda2)$, and can be evaluated by the trace property:
$\RI\lambda2(\phi_l\phi_l^*)= \LI\lambda2(\phi_l\phi_l^*)
=\frac{d(\lambda_1)}{d(\lambda_2)}\LI\lambda1(\phi_l^*\phi_l)=
\frac{d(\lambda_1)}{d(\lambda_2)}$, while the sum over $l$ yields the
multiplicity factor $Z_{\lambda_1,\lambda_2}$. Hence 
$$\sum_{lm,e_1e_2} \overline{\zeta^s_{lm,e_1e_2}}\zeta^n_{lm,e_1e_2} 
  = \left(\sum_{\lambda_1,\lambda_2}
  \frac{d(\lambda_1)d(\lambda_2)Z_{\lambda_1,\lambda_2}}{d(\theta)}\right) 
  \LI\nu1(\phi_s^*\phi_n)=\delta_{sn},$$
and hence $w_1^*w_1=\sum_n W_nW_n^*=1$.
This completes the proof of the Theorem. 
\end{beweis}
\begin{beweis}{ of Proposition 1.5} Left multiplication of $w_1$ with the
induced braiding operator 
$$\eps(\theta,\theta) = \sum_{mlm'l'} (W_{m'}\times W_{l'})
\cir(\eps_1(\lambda_1,\mu_1) \otimes (\eps_2(\lambda_2,\mu_2)^*)\opp)\cir
(W_l\times W_m)^* $$
amounts to a unitary passage from bases $T_e\in Hom(\nu,\lambda\mu)$ to
bases $\eps(\lambda,\mu)T_e\in Hom(\nu,\mu\lambda)$. But by (E3), the
coefficients $\zeta^n_{lm,e_1e_2}$ are invariant under these changes
of bases. Hence $\eps(\theta,\theta)w_1=w_1$. 
\end{beweis}
\section{Joint inclusions and normality}
The main purpose of this section is to introduce and discuss the
notions of ``normality'' and ``essential normality''. These properties
are of interest since the embedding of chiral subtheories into
two-dimensional conformal quantum field theories should always give
rise to essentially normal CTPS's \cite{R1}. We start by
introducing these and related notions in the broader context of
``joint inclusions'' of von Neumann algebras, i.e., triples $(A,B,C)$
such that $A\vee B\subset C$. We first record some more or less
elementary properties of joint inclusions, before we give a simple
characterization of normality in the case of CTPS's in terms of the
coupling matrix. 
\begin{satz}{Definition} 
  Let $\Lambda=(A,B,C)$ be a joint inclusion of 
  von Neumann algebras. We denote by $\Lambda^c:=(B^c,A^c,C)$ the joint 
  inclusion of the relative commutants in $C$. We write $\Lambda_1\subset
  \Lambda_2$ if $C_1=C_2$ and $A_1\subset A_2$, $B_1\subset B_2$, and call
  $\Lambda_2$ {\em intermediate} w.r.t.\ $\Lambda_1$. We call
  $\Lambda$ {\em symmetric} if $\Lambda\subset\Lambda^c$ (i.e., $A$ and $B$
  commute with each other). We call $\Lambda$ {\em normal} if 
  $\Lambda=\Lambda^c$ (i.e., $A$ and $B$ are each other's relative 
  commutants). We call $\Lambda$ {\em essentially normal} if
  $\Lambda^c=\Lambda^{cc}$. 
\end{satz} 
One has the following elementary facts. 
\begin{satz}{Proposition} 1. $\Lambda^c=\Lambda^{ccc}$.\\[2pt]\indent
2. If $\Lambda_1\subset\Lambda_2$ then 
$\Lambda_2^c\subset\Lambda_1^c$.\\[2pt]\indent
3. If $\Lambda$ is symmetric then $\Lambda \subset \Lambda^{cc}
\subset \Lambda^c$.\\[2pt]\indent 
4. $\Lambda$ is essentially normal if and only if $\Lambda$ and $\Lambda^c$
are both symmetric. \\[2pt]\indent
5. Every symmetric $\Lambda$ has a normal intermediate joint
inclusion. \\[2pt]\indent
6. If $\Lambda$ is normal then one has $Z(A)=(A\vee B)^c=Z(B)
\supset Z(C)$, so $A$ and likewise $B$ are factors if and only if 
$A\vee B\subset C$ is irreducible, and in this case $C$ necessarily is 
also a factor.
\end{satz}
\begin{beweis}{}
Assertions 1.--4. are obvious. A normal intermediate joint inclusion
is given by, e.g., $(B^c,B^{cc},C)$. Assertion 6 holds since 
$(A\vee B)^c= A^c \cap B^c$.
\end{beweis}
While these statements are in quite some parallelism to the 
theory of self-adjoint extensions of symmetric unbounded operators,
assertion 5 is a departure from this parallelism, since 
self-adjoint extensions do not always exist for symmetric
operators. The parallelism seems to become closer if one restricts to the
subclass of tensor product subfactors (canonical or not) which are
obviously symmetric joint inclusions. But neither
$((\eins\otimes B)^c,(A\otimes\eins)^c,C)$ nor the joint inclusion
$((\eins\otimes B)^c,(\eins\otimes B)^{cc},C)$ in 
4.2(5) will again be a tensor product subfactor in general.

While we have no general criterium for the existence of normal
intermediate tensor product subfactors in general, the following
proposition gives a simple characterization of normality in the case
of CTPS's, which entails certain constraints on the structure of 
$A_1\otimes B_1\subset C$ for which a normal CTPS $A\otimes B\subset C$ 
can possibly be intermediate. These constraints will apply to the
embeddings of left and right chiral subtheories into two-dimensional
conformal quantum field theories, which by \cite{R1} give rise to
CTPS's whose relative commutants are again tensor product subfactors,
hence symmetric. Thus these local subfactors are essentially normal
CTPS's by Prop.\ 4.2(4), and the normal intermediate subfactor
corresponds to the maximally extended chiral algebras (going along
with permutation modular invariants). We do not evaluate these
constraints here, but it is clear that the total dual canonical
endomorphism must be of the form
$(\bar\iota_A\otimes\bar\iota_B)\cir\theta\cir(\iota_A\otimes\iota_B)
\simeq\bigoplus_\alpha\,(\bar\iota_A\alpha\iota_A)\otimes
(\bar\iota_B\sigma(\alpha)\iota_B)$ where $\theta$ corresponding to
the normal intermediate inclusion is of the special ``permutational''
form (N3) as described in the following proposition. \newpage 
\begin{satz}{Proposition} 
  Let $A \otimes B \subset C$ be a CTPS of type III with coupling 
  matrix $Z$, i.e., the dual canonical endomorphism is of the form
  $$ \theta\simeq \bigoplus_{\alpha\in\Delta_A,\beta\in\Delta_B}
  Z_{\alpha,\beta}\;\alpha\otimes\beta , $$
  where $\Delta_A\ni\id_A$ and $\Delta_B\ni\id_B$ are two sets of mutually 
  inequivalent irreducible endomorphisms in $End\fin(A)$ and
  $End\fin(B)$. Then the following conditions are equivalent. 

  (N1) The joint inclusion $(A\otimes\eins_B,\eins_A\otimes B,C)$ is
  normal, i.e., $A\otimes \eins_B$ and $\eins_A\otimes B$ are each
  other's relative commutants in $C$.  

  (N2) The coupling matrix couples no non-trivial sector of $A$ to the
  trivial sector of $B$, and vice versa, i.e.,
  $$ Z_{\alpha,\id_B}=\delta_{\alpha,\id_A} \quad \hbox{\rm and} \quad
  Z_{\id_A,\beta}=\delta_{\beta,\id_B} . $$

  (N3) The sets $\Delta_A$ and $\Delta_B$ are closed $A$- and $B$-systems, 
  respectively, i.e., they are both closed under conjugation and fusion. 
  There is a bijection $\sigma:\Delta_A \to\Delta_B$ which preserves
  the fusion rules, i.e.,
  $$ \dim Hom(\alpha_1,\alpha_2\alpha_3) = \dim
  Hom(\sigma(\alpha_1),\sigma(\alpha_2)\sigma(\alpha_3)) . $$
  The matrix $Z$ is the permutation matrix for this bijection, i.e.,
  $$ Z_{\alpha,\beta} = \delta_{\sigma(\alpha),\beta} . $$
\end{satz}
The proof is published in \cite[Lemma 3.4 and Thm.\ 3.6]{R1}. 

Making contact with the new CTPS's in Thm.\ 1.4, we first point
out that in general they will {\em not} be normal, since among the coupling
matrices constructed in \cite{BEK} there are those which are not of
the form (N2,N3). 

The most simple case $N_1=N_2=M$ hence $Z=\eins$ is known for a while
\cite{LR}, and clearly {\em is} normal by Prop.\ 4.3. Specifically, it
describes the ``diagonal'' extension of the chiral observables by local
two-dimensional observables carrying opposite chiral charges. We
conclude from Prop.\ 4.3 that the left and right chiral observables
are each other's relative commutants within this two-dimensional
theory, and the same holds whenever the coupling matrix in Cor.\ 1.6
satisfies condition (N2,N3), i.e., describes a permutation modular
invariant (cf.\ \cite{R1}).

In the abstract mathematical setting, the subfactors with $Z=\eins$
constructed in \cite{LR} were recognized \cite{M} (up to some trivial
tensoring with a type III factor) as the type II asymptotic subfactor
\cite{O} associated with $\sigma(M)\subset M$ where 
$\sigma\equiv\bigoplus_{\lambda\in\Delta}\lambda$. As the 
asymptotic subfactor $M \vee M^c\subset M_\infty$ associated with a  
fixed point inclusion $M^G\subset M$ for an outer action of a group
$G$, provides the same category of $M_\infty$-$M_\infty$ bimodules as
a fixed point inclusion for an outer action of the quantum double
$D(G)$ on $M_\infty$, general asymptotic subfactors in turn are
considered \cite{O,EK} as generalized quantum doubles. 

General asymptotic subfactors are CTPS's, i.e., $M\vee M^c \simeq
M\otimes M^c$ are in a tensor product position within $M_\infty$, and
every irreducible $M\vee M^c$-$M\vee M^c$ bimodule associated with the
asymptotic subfactor respects the tensor product \cite{O}. They are
normal, i.e., $M$ and $M^c$ are each other's relative commutant in 
$M_\infty$. Moreover, the system of $M_\infty$-$M_\infty$ bimodules
associated with an asymptotic subfactor has a non-degenerate braiding
\cite{O,I}. We do not know at present whether the new CTPS's always
share this braiding property, which ought to be tested with methods as
in \cite{I}. 

\section{Conclusion}
We have shown the existence of a class of new subfactors associated
with extensions of closed systems of sectors. The proof proceeds by
establishing the corresponding Q-systems in terms of certain matrix
elements for the transition between two extensions. The new subfactors
are canonical tensor product subfactors and include the asymptotic
subfactors.  

Among the new subfactors, there are the local subfactors of two-dimensional 
conformal quantum field theory associated with certain modular invariants, 
thereby establishing the expected existence of these theories. We also
gave a characterization of normality of CTPS's which corresponds to
the maximal subtheories of chiral observables in these models.

\begin{acknowledgement}
I am deeply indebted to Y. Kawahigashi, M. Izumi, T. Matsui, and
I. Ojima who made possible my visit to Japan during summer 1999 where
the present work was completed. I want to thank all of them as well as
H. Kosaki, Y. Watatani, and T. Masuda for discussions, and for their
hospitality extended to me at the Department of Mathematical Sciences,
University of Tokyo, the Graduate School of Mathematics, Kyushu
University, and the Research Institute for Mathematical Sciences,
Kyoto University. I also thank H. Kurose for giving me the opportunity
to present these results at the workshop ``Advances in Operator
Algebras'' \cite{R2} held at RIMS, Kyoto. Financial support by a
Grant-in-Aid for Scientific Research from the Ministry of Education
(Japan) is gratefully acknowledged. Finally, I thank J. B\"ockenhauer for
sending me a preliminary manuscript on related issues from a
complementary perspective \cite{BEp}. 
\end{acknowledgement}

\small \addtolength{\baselineskip}{-.5pt}

\end{document}